# KOLMOGOROV EQUATION ASSOCIATED TO THE STOCHASTIC REFLECTION PROBLEM ON A SMOOTH CONVEX SET OF A HILBERT SPACE

By Viorel Barbu,[1] Giuseppe Da Prato[2] and Luciano Tubaro[2]

*University Al. I. Cuza, Scuola Normale Superiore and University of Trento*

We consider the stochastic reflection problem associated with a self-adjoint operator $A$ and a cylindrical Wiener process on a convex set $K$ with nonempty interior and regular boundary $\Sigma$ in a Hilbert space $H$. We prove the existence and uniqueness of a smooth solution for the corresponding elliptic infinite-dimensional Kolmogorov equation with Neumann boundary condition on $\Sigma$.

**1. Introduction.** Let us consider a stochastic differential inclusion in a Hilbert space $H$,

$$(1.1) \qquad \begin{cases} dX(t) + (AX(t) + N_K(X(t)))\,dt \ni dW(t), \\ X(0) = x. \end{cases}$$

Here $A : D(A) \subset H \to H$ is a self-adjoint operator, $K = \{x \in H : g(x) \le 1\}$, where $g : H \to \mathbb{R}$ is convex and of class $C^\infty$, $N_K(x)$ is the normal cone to $K$ at $x$ and $W(t)$ is a cylindrical Wiener process in $H$ (see Hypothesis 1.1 for more precise assumptions). Obviously the expression in (1.1) is formal and its precise meaning should be defined.

When $H$ is finite-dimensional a solution to (1.1) is a pair of continuous adapted processes $(X, \eta)$ such that $X$ is $K$-valued, $\eta$ is of bounded variation with $d\eta$ concentrated on the set of times where $X(t) \in \Sigma$ (the boundary of $K$) and

$$X(t) + \int_0^t AX(s)\,ds + \eta(t) = x + W(t), \qquad t \ge 0, \ \mathbb{P}\text{-a.s.},$$

Received May 2008; revised September 2008.

[1]Supported by Grant PN-II ID 404 (2007–2010) of Romanian Minister of Research.

[2]Supported in part by the Italian National Project MURST "Equazioni di Kolmogorov."

*AMS 2000 subject classifications.* 60J60, 47D07, 15A63, 31C25.

*Key words and phrases.* Reflected process, convex sets, Dirichlet forms, Kolmogorov operators, Gaussian measures, infinite-dimensional Neumann problem.

$$\int_0^T (d\eta(t), X(t) - z(t)) \geq 0, \qquad \mathbb{P}\text{-a.s.,}$$

for all $z \in C([0,T]; K)$. The existence and uniqueness of a solution $(X, \eta)$ to latter equation was first proven by Cépa in [5]. (See also [3] for a slightly different formulation.)

Therefore, under the assumptions of [3] or [5], one can construct a transition semigroup in $C(K)$ by the usual formula

$$P_t\varphi(x) = \mathbb{E}[\varphi(X(t,x))], \qquad t \geq 0, \ \varphi \in C(K).$$

The infinitesimal generator $L$ of $P_t$ is the Kolmogorov operator

$$L\varphi = \tfrac{1}{2}\Delta\varphi + \langle Ax, D\varphi \rangle$$

equipped with a Neumann condition at the boundary $\Sigma$ of $K$. (See, e.g., [3], where the more general case of oblique derivative boundary conditions were also considered.)

Let us go now to the infinite-dimensional situation. In this context (1.1) was first studied by Nualart and Pardoux [18], when $H = L^2(0,1)$, $A$ is the Laplace operator with Dirichlet or Neumann boundary conditions and $K$ is the convex set of all nonnegative functions of $L^2(0,1)$; see also [13].

The Kolmogorov operator in this situation was described by Zambotti [21], in the space $L^2(H, \nu)$ where $\nu$ is the law of the $3D$-Bessel Bridge which coincides with the unique invariant measure of (1.1). Zambotti was able to show that the Dirichlet form

$$a(u, v) = \int_K \langle Du, Dv \rangle \, d\nu$$

is closable by proving a suitable integration by parts formula and to construct the corresponding Markov semigroup.

Except the situation mentioned above, no existence and uniqueness results for (1.1) are known for the infinite-dimensional equation (1.1). Also it was so far open the characterization of the of the domain of the corresponding Kolmogorov operator.

In this paper we shall consider a regular convex set $K$ with nonempty interior and, though this does not cover the case considered by [21], we are able, however, to get sharp informations on the Kolmogorov generator for a quite general class of convex sets $K$. In this way, though we are not able to approach directly the stochastic variational problem (1.1), we can instead find a regular solution of the corresponding infinite-dimensional Kolmogorov equation equipped with the Neumann boundary condition,

(1.2) $$\begin{cases} \lambda\varphi - \tfrac{1}{2}\mathrm{Tr}[D^2\varphi] - \langle x, AD\varphi \rangle = f, & x \in K, \\ \langle D\varphi, N_K(x) \rangle = 0, & \forall x \in \Sigma, \end{cases}$$

where $\lambda > 0$ and $f \in L^2(K, \nu)$.

In this way we obtain a Markov semigroup $P_t$ which by the results of [16] provides a process corresponding to a martingale solution of (1.1) (see also the forthcoming paper [1]).

A basic tool we are using is a co-area formula from Malliavin; see [17] valid for $g$ of class $C^\infty$. Moreover, in the Appendix we present a direct proof of this formula when $g$ is $C^2$ and fulfills some additional conditions which are covered in several situations, for instance when $K$ is a ball; in that case the co-area formula was proved (1979) by Hertle [14].

Let us explain the content of this paper. As we said, we take a convex set of the form $K = \{x \in H : g(x) \le 1\}$ where $g : H \to \mathbb{R}$ is of class $C^\infty$ and with second order derivative $D^2 g$ positive definite. Then we consider the probability measure $\nu$ given for any Borel set $I$ of $K$ by

$$\nu(I) = \frac{\mu(I)}{\mu(K)},$$

where $\mu$ is the Gaussian measure (corresponding to the linear problem without reflection) of mean 0 and covariance $Q = \frac{1}{2} A^{-1}$.

In Section 2, by exploiting a basic infinite-dimensional co-area formula, see [17], we are able to prove an integration by parts formula for $\nu$. This allows us to show in Section 3 that the Dirichlet form

$$a(u, v) = \int_K \langle Du, Dv \rangle \, d\nu$$

is closable (see also [1] for a different approach). In this way, by the usual variational theory, we can define its generator $N$ and construct the corresponding Markov transition semigroup $P_t$, which is reversible since $N$ is self adjoint.

In Section 4 we study the Kolmogorov equation (1.2) by the classical method of penalization

$$\lambda \varphi_\varepsilon - \frac{1}{2} \operatorname{Tr}[D^2 \varphi_\varepsilon] + \langle x, AD\varphi_\varepsilon \rangle + \frac{1}{\varepsilon} \langle x - \Pi_K(x), D\varphi_\varepsilon \rangle = f, \qquad x \in H,$$
(1.3)

where $\Pi_K(x)$ is the projection of $x$ on $K$. We show that $\{\varphi_\varepsilon\}$ strongly converges to the solution $\varphi = (\lambda I - N)^{-1} f$ of (1.2) and that

$$(1.4) \qquad D(N) \subset \left\{ \varphi \in W^{2,2}(K, \nu) : \int_K |A^{1/2} D\varphi|^2 \, d\nu < +\infty \right.$$
$$\left. \text{and } \langle D\varphi, N_K(x) \rangle = 0 \text{ on } \Sigma \right\}.$$

These results seem to be new in infinite dimensions; see [2, 3, 7] for the finite-dimensional case.

Finally, Section 5 is devoted to equations of the form

(1.5)    $\begin{cases} dX(t) + (AX(t) + F(X(t)) + N_K(X(t))) \, dt \ni dW(t), \\ X(0) = x, \end{cases}$

where $F \colon H \to H$ is a nonlinear perturbation of $A$.

In Section 5.1 we assume that $F = DV$ where $V \colon H \to \mathbb{R}$ is a regular potential. This case is an easy generalization of the previous one (i.e., when $F = 0$), namely measure $\nu$ is replaced by the following one:

$$\zeta(dx) = \frac{e^{-2V(x)}}{\int_K e^{-2V(y)} \nu(dy)} \nu(dx).$$

This extension is briefly described in that section.

In Section 5.2 the case of a bounded Borel function $F$, not necessarily of potential type, is considered. Here we can solve the Kolmogorov equation

(1.6)    $\begin{cases} \lambda \varphi - \frac{1}{2} \operatorname{Tr}[D^2 \varphi] + \langle x, AD\varphi \rangle - \langle F(x), D\varphi \rangle = f, & x \in K, \\ \langle D\varphi, N_K(x) \rangle = 0, & \forall x \in \Sigma \end{cases}$

by a straightforward perturbation argument, taking avantage of the inclusion (1.4). In this way we obtain a solution $\varphi \in D(N)$ of (1.6) only for $\lambda$ sufficiently large. Also, obviously, measure $\nu$ is not invariant for the corresponding semigroup $Q_t$. However, using the fact that operator $Q_t$ is compact in $L^2(K, \nu)$, we can show the existence of an invariant measure $\zeta$ for $Q_t$ so that the extension of $Q_t$ to $L^1(K, \zeta)$ is the natural transition semigroup associated with (1.5). Notice, however, that this semigroup is not reversible (when $F$ is not of potential type).

We conclude this section by precising assumptions and notation which will be used throughout in what follows.

*Assumptions.*  We are given a real separable Hilbert space $H$ (with scalar product $\langle \cdot, \cdot \rangle$ and norm denoted by $|\cdot|$). Concerning $A$, $K$ and $W$ we shall assume that:

HYPOTHESIS 1.1.  (i) $A \colon D(A) \subset H \to H$ is a linear self-adjoint operator on $H$ such that $\langle Ax, x \rangle \geq \delta |x|^2, \forall x \in D(A)$ for some $\delta > 0$. Moreover, $A^{-1}$ is of trace class.

(ii) There exists a convex $C^\infty$ function $g \colon H \to \mathbb{R}$ with $D^2 g$ positively defined, that is, $\langle D^2 g(x)h, h \rangle \geq \gamma |h|^2$, $\forall h \in H$ where $\gamma > 0$, such that

$$K = \{x \in H : g(x) \leq 1\}, \qquad \Sigma = \{x \in H : g(x) = 1\}.$$

(iii) $W$ is a cylindrical Wiener process on $H$ of the form

$$W(t) = \sum_{k=1}^{\infty} \beta_k(t) e_k, \qquad t \geq 0,$$

where $\{\beta_k\}$ is a sequence of mutually independent real Brownian motions on a filtered probability spaces $(\Omega, \mathcal{F}, \{\mathcal{F}_t\}_{t\geq 0}, \mathbb{P})$ (see, e.g., [8]) and $\{e_k\}$ is an orthonormal basis in $H$ which will be taken as a system of eigen-functions for $A$ for simplicity, that is,

$$Ae_k = \alpha_k e_k \qquad \forall k \in \mathbb{N},$$

where $\alpha_k \geq \delta$.

We notice that the interior $\overset{\circ}{K}$ is nonempty since $D^2 g$ is positive definite.

*Notation.* We denote by $\mathcal{B}(H)$ [resp. $\mathcal{B}(K)$] the $\sigma$-field of all Borel subsets of $H$ (resp. $K$) and by $\mathcal{P}(H)$ [resp. $\mathcal{P}(K)$] the set of all probability measures on $(H, \mathcal{B}(H))$ [resp. $(K, \mathcal{B}(K))$].

Everywhere in the following $D\varphi$ is the derivative of a function $\varphi \colon H \to \mathbb{R}$. By $D^2\varphi \colon H \to L(H, H)$ we shall denote the second derivative of $\varphi$. We shall denote also by $C_b(H)$ and $C_b^k(H), k \in \mathbb{N}$, the spaces of all continuous and bounded functions on $H$ and, respectively, of $k$-times differentiable functions with continuous and bounded derivatives. The space $C^k(K), k \in \mathbb{N}$, is defined as the space of restrictions of functions of $C_b^k(H)$ to the subset $K$.

The boundary of $K$ will be denoted by $\Sigma$. $N_K(x)$ is the normal cone to $K$ at $x$, that is,

$$N_K(x) = \{z \in H : \langle z, y - x \rangle \leq 0, \forall y \in K\}.$$

Moreover, we shall denote by $d_K(x)$ the distance of $x$ from $K$ and by $I_K$ the indicator function of $K$,

$$I_K(x) = \begin{cases} 0, & \text{if } x \in K, \\ +\infty, & \text{if } x \notin K. \end{cases}$$

For any $\varepsilon > 0$, $U_\varepsilon$ will represent the Moreau approximation of $I_K$ given by

$$U_\varepsilon(x) = \inf\left\{I_K(y) + \frac{1}{2\varepsilon}|x - y|^2, y \in H\right\} = \frac{1}{2\varepsilon}d_K(x)^2, \qquad x \in H.$$

It is well known that

$$DU_\varepsilon(x) = \frac{1}{\varepsilon}(x - \Pi_K(x)), \qquad x \in H, \varepsilon > 0,$$

where $\Pi_K(x)$ is the projection of $x$ over $K$. In particular, we have

$$(1.7) \qquad D(d_K^2(x)) = x - \Pi_K(x) \qquad \forall x \in K^c,$$

($K^c$ is the complement of $K$) which implies

$$(1.8) \qquad Dd_K(x) = \frac{x - \Pi_K(x)}{d_K(x)} \qquad \forall x \in K^c.$$

We denote by $\mathbf{n}(\Pi_K(x))$ the exterior normal at $\Pi_K(x)$,

$$\mathbf{n}(\Pi_K(x)) = \frac{x - \Pi_K(x)}{d_K(x)} \qquad \forall x \in K^c.$$

From (1.8) we deduce that

$$(1.9) \quad D(x - \Pi_K(x)) = Dd_K(x) \otimes Dd_K(x) + d_K(x)D^2 d_K(x) \qquad \forall x \in K^c.$$

Finally, $\mu$ will represent the Gaussian measure in $H$ with mean 0 and covariance operator

$$Q := \tfrac{1}{2}A^{-1}.$$

Since $A$ is strictly positive $\mu$ is nondegenerate and full. We set

$$\lambda_k = \frac{1}{2\alpha_k} \qquad \forall k \in \mathbb{N},$$

so that

$$Qe_k = \lambda_k e_k \qquad \forall k \in \mathbb{N}.$$

We denote by $\mathcal{E}_A(H)$ the space of all real and imaginary parts of exponential functions $e^{i\langle h, x \rangle}, h \in D(A)$. Then the operator $D \colon \mathcal{E}_A(H) \subset L^2(H, \mu) \to L^2(H, \mu; H)$ is closable in $L^2(H, \mu)$ and the domain of its closure is denoted by $W^{1,2}(H, \mu)$ (the Sobolev space).

The following integration by parts formula for the measure $\mu$ is well known (see, e.g., [9]). For any $\varphi, \psi \in W^{1,2}(H, \mu)$ and $z \in H$,

$$(1.10) \quad \int_H \langle D\varphi, Q^{1/2}z \rangle \psi \, d\mu = -\int_H \langle D\psi, Q^{1/2}z \rangle \varphi \, d\mu + \int_H W_z \varphi \psi \, d\mu,$$

where $W_z$ represents the *white noise* function,

$$(1.11) \quad W_z(x) = \sum_{k=1}^{\infty} \frac{1}{\sqrt{\lambda_k}} \langle x, e_k \rangle \langle z, e_k \rangle \qquad \forall z \text{ and } \mu\text{-a.e. } x \in H.$$

We recall that $W_z$ is a Gaussian random variable in $L^2(H, \mu)$ with mean 0 and covariance $|z|^2$.

**2. The measure $\boldsymbol{\mu}$ conditioned to $\boldsymbol{K}$.** We denote by $\nu$ the Gaussian measure $\mu$ conditioned to $K$, that is,

$$\nu(I) = \frac{\mu(K \cap I)}{\mu(K)} \qquad \forall I \in \mathcal{B}(H).$$

Since $\mu$ is full and $\overset{\circ}{K}$ is nonempty, this definition is meaningful. We notice that, thanks to Hypothesis 1.1(ii) the surface measure $\mu_\Sigma$ is well defined (see [17]).

We want now to prove an integration by parts formula with respect to measure $\nu$ which generalizes (1.10). For this it is convenient to introduce a sequence of approximating measures $\{\nu_\varepsilon\}_{\varepsilon>0}$ defined by,

$$(2.1) \qquad \nu_\varepsilon(dx) = \rho_\varepsilon(x)\mu(dx), \qquad x \in H,$$

where,

$$(2.2) \qquad \rho_\varepsilon(x) = Z_\varepsilon^{-1} e^{-1/\varepsilon\, d_K^2(x)}$$

and

$$(2.3) \qquad Z_\varepsilon = \int_H e^{-1/\varepsilon\, d_K^2(y)} \mu(dy).$$

Notice that, by the dominated convergence theorem,

$$(2.4) \qquad \lim_{\varepsilon \to 0} Z_\varepsilon = Z_0 = \mu(K),$$

whereas

$$(2.5) \qquad \lim_{\varepsilon \to 0} \rho_\varepsilon(x) = \begin{cases} 1, & \text{if } x \in K, \\ 0, & \text{if } x \notin K. \end{cases}$$

So, we have

$$(2.6) \qquad \lim_{\varepsilon \to 0} \nu_\varepsilon = \nu \qquad \text{weakly in } \mathcal{P}(H).$$

Moreover

$$(2.7) \qquad D\rho_\varepsilon(x) = -\frac{2}{\varepsilon} \rho_\varepsilon(x)(x - \Pi_K(x)),$$

so that $\rho_\varepsilon \in W^{1,2}(H,\mu)$.

We shall denote by $L^2(K,\nu)$ the space of all $\nu$-square-integrable functions on $K$ with the scalar product

$$\langle u, v \rangle_{L^2(K,\nu)} = \int_K u(x)v(x)\nu(dx)$$

and the norm $|u|^2_{L^2(K,\nu)} = \langle u, u \rangle_{L^2(K,\nu)}$.

2.1. *The integration by parts formula.* Here we are going to derive from (1.10), an integration by parts formula for the measure $\nu_\varepsilon$. Let $\varphi \in C_b^1(H), z \in H$, then, since $\rho_\varepsilon \in W^{1,2}(H,\mu)$, we find from (1.10) that

$$\int_H \langle D\varphi, Q^{1/2}z \rangle \, d\nu_\varepsilon = \int_H \langle D\varphi, Q^{1/2}z \rangle \rho_\varepsilon \, d\mu$$

$$= -\int_H \varphi \langle D \log \rho_\varepsilon, Q^{1/2}z \rangle \, d\nu_\varepsilon + \int_H W_z \varphi \, d\nu_\varepsilon.$$

Since,

$$D \log \rho_\varepsilon(x) = -\frac{1}{\varepsilon}(x - \Pi_K x),$$

we find the formula

$$
\begin{aligned}
(2.8) \quad \int_H \langle D\varphi, Q^{1/2}z \rangle \nu_\varepsilon(dx) &= \frac{1}{\varepsilon} \int_H \varphi(x) \langle x - \Pi_K(x), Q^{1/2}z \rangle \nu_\varepsilon(dx) \\
&\quad + \int_H W_z(x)\varphi(x)\nu_\varepsilon(dx).
\end{aligned}
$$

LEMMA 2.1. *Let* $\varphi \in C_b^1(H), z \in H.$ *Then there exists the limit*

$$
\begin{aligned}
(2.9) \quad \lim_{\varepsilon \to 0} J_\varepsilon^z(\varphi) &:= \lim_{\varepsilon \to 0} \frac{1}{\varepsilon} \int_H \varphi(x) \langle x - \Pi_K x, Q^{1/2}z \rangle \nu_\varepsilon(dx) \\
&= \int_\Sigma \varphi(y) \langle \mathbf{n}(y), Q^{1/2}z \rangle \mu_\Sigma(dy),
\end{aligned}
$$

*where* $\mathbf{n}(y) = \nabla g(y)/|\nabla g(y)|$ *is the exterior normal to* $\Sigma$ *at* $y$ *and* $\mu_\Sigma$ *is the surface measure on* $\Sigma$ *induced by* $\mu$ *(see [17]).*

PROOF. First we notice that

$$J_\varepsilon^z(\varphi) = \frac{1}{\varepsilon Z_\varepsilon} \int_{\{d_K(x)>0\}} \varphi(x) d_K(x) \langle \mathbf{n}(\Pi_K(x)), Q^{1/2}z \rangle e^{-d_K^2(x)/\varepsilon} \mu(dx).$$

By the co-area formula (see [17], page 140) (see also Theorem A.5 below) we have

$$(2.10) \qquad \int_H f\mu(dx) = \int_0^\infty \left[ \int_{\Sigma_r} f(y)\mu_{\Sigma_r}(dy) \right] dr.$$

Notice that the surface measure is defined for all $r \geq 0$ taking into account ([17], Theorem 6.2, Chapter V); moreover ([17], Theorem 1.1, Corollary 6.3.2, Chapter V), give the continuity property in Theorem 6.3.1 of [17], Chapter V. Setting in (2.10)

$$f = (1 - \mathbb{1}_K)\varphi(x) \, d_K(x) \langle \mathbf{n}(\Pi_K(x)), Q^{1/2}z \rangle e^{-d_K^2(x)/\varepsilon},$$

we get

$$J_\varepsilon^z(\varphi) = \frac{1}{\varepsilon Z_\varepsilon} \int_0^{+\infty} \xi e^{-\xi^2/\varepsilon} \, d\xi \int_{\Sigma_{\xi+1}} \varphi(y) \langle \mathbf{n}(\Pi_K(x)), Q^{1/2}z \rangle \mu_{\Sigma_\xi}(dy).$$

Hence, setting $\xi = \sqrt{\varepsilon}s$, yields

$$J_\varepsilon^z(\varphi) = \frac{1}{Z_\varepsilon} \int_0^\infty s e^{-s^2} \, ds \int_{\Sigma_{\sqrt{\varepsilon}s+1}} \varphi(y) \langle \mathbf{n}(\Pi_K(y)), Q^{1/2}z \rangle \mu_{\Sigma_{\sqrt{\varepsilon}s}}(dy).$$

So (2.9) follows. $\quad\square$

We are now in position to prove the announced integration by parts formula.

THEOREM 2.2.    *Let $\varphi \in C_b^1(H), z \in H$. Then for any $z \in H$ we have*

$$\begin{aligned}
(2.11) \quad \int_K \langle D\varphi(x), Q^{1/2}z\rangle \nu(dx) &= \frac{1}{2\mu(K)} \int_\Sigma \varphi(y)\langle \mathbf{n}(y), Q^{1/2}z\rangle \mu_\Sigma(dy) \\
&\quad + \int_K W_z(x)\varphi(x)\nu(dx).
\end{aligned}$$

PROOF.    The conclusion of the theorem follows letting $\varepsilon \to 0$ in (2.8) and taking into account Lemma 2.1. $\quad\square$

2.2. *The Sobolev space $W^{1,2}(K,\nu)$.*    We shall define space $W^{1,2}(K,\nu)$ by proving, as it is usual, closability of the gradient. For this we need a lemma.

LEMMA 2.3.    *The space*

$$C_0^1(K) := \{\varphi \in C^1(K) : \varphi = 0 \ on \ \Sigma\}$$

*is dense in $L^2(K,\nu)$.*

PROOF.    It is enough to show that if $\varphi \in C^1(K)$ then there exists a sequence $\{\varphi_\alpha\} \subset C_0^1(K)$ such that

$$(2.12) \qquad\qquad \lim_{\alpha \to 0} \varphi_\alpha = \varphi \qquad \text{in } L^2(K,\nu).$$

Let $\{\chi_\alpha\}_{\alpha \in (0,1)} \subset C^1(\mathbb{R})$ be a sequence such that,

$$\chi_\alpha(r) = \begin{cases} 1, & \text{for } r \in [0, 1-\alpha], \\ 0, & \text{for } r \geq 1. \end{cases}$$

Setting now

$$\varphi_\alpha(x) = \chi_\alpha(g(x))\varphi(x) \qquad \forall \alpha \in (0,1),$$

we see that $\{\varphi_\alpha\}_{\alpha \in (0,1)} \subset C_0^1(K)$ and (2.12) follows from the dominated convergence theorem. $\quad\square$

PROPOSITION 2.4.    *The mapping*

$$D : C^1(K) \subset L^2(K,\nu) \to L^2(K,\nu;H), \qquad \varphi \to D\varphi,$$

*is closable.*

PROOF. Let $(\varphi_n) \subset C^1(K)$ be such that

$$\varphi_n \to 0 \quad \text{in } L^2(K, \nu), \qquad D\varphi_n \to F \quad \text{in } L^2(K, \nu; H)$$

as $n \to \infty$. We have to show that $F = 0$. Let $\psi \in C_0^1(K)$ and $z \in H$. Then by (2.11) with $\varphi_n \psi$ replacing $\varphi$ (see Theorem 2.2) we have that

$$
\begin{aligned}
& \int_K \langle D\varphi_n(x), Q^{1/2}z\rangle \psi(x)\nu(dx) \\
(2.13) \quad & = -\int_K \langle D\psi(x), Q^{1/2}z\rangle \varphi_n(x)\nu(dx) \\
& \quad + \frac{1}{2\mu(K)} \int_\Sigma \varphi_n(y)\psi(y)\langle \mathbf{n}(y), Q^{1/2}z\rangle \mu_\Sigma(dy) \\
& \quad + \int_K W_z(x)\varphi_n(x)\psi(x)\nu(dx) \\
& = -\int_K \langle D\psi(x), Q^{1/2}z\rangle \varphi_n(x)\nu(dx) + \int_K W_z(x)\varphi_n(x)\psi(x)\nu(dx),
\end{aligned}
$$

since $\psi$ vanishes on $\Sigma$. Letting $n \to \infty$ we find that

$$\int_H \langle F(x), Q^{1/2}z\rangle \psi(x)\mu(dx) = 0.$$

This implies $F = 0$ in view of the arbitrariness of $\psi$ and $z$ [recall Lemma 2.3 and that $Q^{1/2}(H)$ is dense in $H$]. $\square$

We shall still denote by $D$ the closure of $D$ and by $W^{1,2}(K, \nu)$ its domain of definition. $W^{1,2}(K, \nu)$ is a Hilbert space with the scalar product

$$\langle \varphi, \psi\rangle_{W^{1,2}(K,\nu)} = \int_K [\varphi\psi + \langle D\varphi, D\psi\rangle]\, d\nu.$$

2.3. *The trace of a function of $W^{1,2}(K, \nu)$.* In order to define the trace of a function $\varphi \in W^{1,2}(K, \nu)$ we need a technical lemma.

LEMMA 2.5. *Assume that $\varphi \in C_b^1(H)$. Then the following estimate holds,*

$$
\begin{aligned}
(2.14) \quad & \int_\Sigma |Q^{1/2}\mathbf{n}(y)|^2 \varphi^2(y)\mu_\Sigma(dy) \\
& \leq C\left( \int_K \varphi^2(x)\nu(dx) + \int_K |D\varphi(x)|^2\nu(dx) \right),
\end{aligned}
$$

*where $C$ is a suitable constant.*

Proof. Here we follow [12]. Let $\varphi \in C^1(K)$. Set $F(x) = Dg(x)$. In particular $F(x) = |Dg(x)|\mathbf{n}(x)$ for $x \in \Sigma$. Then, replacing in (2.11) $\varphi$ with $\lambda_k F_k \varphi^2$ and $z$ with $e_k$, one gets

$$\int_K \lambda_k D_k F_k \varphi^2 \, d\nu + 2\lambda_k \int_K F_k \varphi D_k \varphi \, d\nu$$

$$= \frac{1}{2\mu(K)} \int_\Sigma \lambda_k |Dg(y)| \langle \mathbf{n}(y), e_k \rangle^2 \varphi^2(y) \mu_\Sigma(dy) + \int_K x_k F_k \varphi^2 \nu(dx).$$

It follows that

$$\frac{1}{2\mu(K)} \int_\Sigma \lambda_k |Dg(y)| \langle \mathbf{n}(y), e_k \rangle^2 \varphi^2(y) \mu_\Sigma(dy)$$

$$\leq \int_K \lambda_k D_k^2 g \varphi^2 \, d\nu + \frac{1}{2} \int_K F_k^2 \varphi^2 \nu(dx) + \frac{1}{2} \lambda_k^2 \int_K |D_k \varphi|^2 \nu(dx)$$

$$- \int_K x_k F_k \varphi^2 \nu(dx).$$

Now the conclusion follows summing up over $k$, since $|Dg|$ is bounded below on $\Sigma$. $\square$

Now we can define the *trace* $\gamma(\varphi)$ on $\Sigma$ of a function $\varphi \in W^{1,2}(K, \nu)$. Let us consider a sequence $\{\varphi_n\} \subset C^1(K)$ strongly convergent to $\varphi$ in $W^{1,2}(K, \nu)$. Then by (2.14) it follows that the sequence $\{|Q^{1/2}\mathbf{n}(y)|(\varphi_n)_\Sigma\}$ is convergent in $L^2(\Sigma, \mu_\Sigma)$ to some element $\tilde{\gamma}(\varphi) \in L^2(\Sigma, \mu_\Sigma)$. Then we set

$$\gamma(\varphi)(y) = \frac{1}{|Q^{1/2}\mathbf{n}(y)|} \tilde{\gamma}(\varphi)(y), \qquad \mu_\Sigma\text{-a.s.}$$

By inequality (2.14) it follows that this definition is consistent, that is, is independent of the sequence $\{\varphi_n\}$ and the map $\varphi \to |Q^{1/2}\mathbf{n}(y)|\gamma(\varphi)$ is continuous from $W^{1,2}(K, \nu) \to L^2(\Sigma, \mu_\Sigma)$. Notice also that though $|Q^{1/2}\mathbf{n}(y)| > 0$ for all $y \in \Sigma$ it is not however bounded from below in infinite dimensions. Now the following result is an immediate consequence of Lemma 2.5 and the density of $C_b^1(H)$ in $W^{1,2}(K, \nu)$.

Proposition 2.6. *Assume that $\varphi \in W^{1,2}(K, \nu)$. Then:*

(i) $|Q^{1/2}\mathbf{n}(y)|\gamma(\varphi) \in L^2(\Sigma, \mu_\Sigma)$,

(ii) *the following estimate holds,*

$$(2.15) \quad \int_\Sigma |Q^{1/2}\mathbf{n}(y)|^2 \varphi^2(y) \mu_\Sigma(dy)$$

$$\leq C \left( \int_K \varphi^2(x) \nu(dx) + \int_K |D\varphi(x)|^2 \nu(dx) \right).$$

We notice that if $H$ is finite-dimensional and $Q = I$ formula (2.15) reduces to a classical result since $|Q^{1/2}\mathbf{n}(y)| = 1$ on $\Sigma$.

2.4. *Compactness of embedding* $W^{1,2}(K,\nu) \subset L^2(K,\nu)$. We first show the log-Sobolev estimate for $\nu$.

PROPOSITION 2.7. *For all* $\varphi \in W^{1,2}(H,\nu)$ *we have*

$$(2.16) \quad \int_K \varphi^2 \log(\varphi^2) \, d\nu \le \frac{1}{\lambda_1} \int_H |D\varphi|^2 \, d\nu + \|\varphi\|^2_{L^2(H,\nu)} \log(\|\varphi\|^2_{L^2(H,\nu)}).$$

PROOF. It is enough to show (2.16) for $\varphi \in C^1(H)$. By [6] (see also [9] and [10]) we know that the log-Sobolev estimate holds for the measure $\nu_\varepsilon$,

$$(2.17) \quad \int_H \varphi^2 \log(\varphi^2) \, d\nu_\varepsilon \le \frac{1}{\lambda_1} \int_H |D\varphi|^2 \, d\nu_\varepsilon + \|\varphi\|^2_{L^2(H,\nu_\varepsilon)} \log(\|\varphi\|^2_{L^2(H,\nu_\varepsilon)}).$$

Now the conclusion follows by (2.6) letting $\varepsilon$ tend to 0.  □

We can now prove the following result.

PROPOSITION 2.8. *The embedding* $W^{1,2}(K,\nu) \subset L^2(K,\nu)$ *is compact.*

PROOF. Let $\{\varphi_n\}$ be a sequence in $W^{1,2}(K,\nu)$ such that

$$(2.18) \qquad\qquad \int_K (\varphi_n^2 + |D\varphi_n|^2) \, d\nu \le C.$$

We have to show that there exists a subsequence of $\{\varphi_n\}$ convergent in $L^2(K,\nu)$. For this we proceed as in [6] noticing that, thanks to the log-Sobolev inequality (2.16), $\{\varphi_n\}$ is uniformly integrable and so, it is enough to find a subsequence of $\{\varphi_n\}$ pointwise convergent to an element of $L^2(K,\nu)$. Let $\{\chi_\alpha\}_{\alpha \in (0,1)} \subset C^1(\mathbb{R})$ be such that,

(i) we have

$$\chi_\alpha(r) = \begin{cases} 1, & \text{for } r \in [0, 1-2\alpha], \\ 0, & \text{for } r \ge 1-\alpha. \end{cases}$$

(ii) $|\chi'_\alpha(r)| \le \frac{2}{\alpha}, \forall \alpha > 0.$

Set now

$$\varphi_n^\alpha(x) = \chi_\alpha(g(x))\varphi_n(x) \qquad \forall \alpha \in (0, 1/2).$$

We claim that for each $\alpha \in (0, 1/2)$ the sequence $\{\varphi_n^\alpha\}_{n \in \mathbb{N}}$ is bounded in $W^{1,2}(H,\mu)$. We have in fact

$$\int_H |\varphi_n^\alpha|^2 \, d\mu = \int_H |\varphi_n^\alpha|^2 \, d\nu \le C$$

and, since

$$D\varphi_n^\alpha(x) = \chi_\alpha(g(x))D\varphi_n(x) + \chi'_\alpha(g(x))\varphi_n(x)Dg(x),$$

we have

$$|D\varphi_n^\alpha(x)| \le |D\varphi_n(x)| + \frac{2}{\alpha}|Dg|_\infty|\varphi_n(x)|.$$

Therefore, there is a positive constant $C_\alpha'$ such that

$$\int_H |D\varphi_n^\alpha|^2 \, d\mu \le C_\alpha'.$$

Recalling that the embedding $W^{1,2}(H,\mu) \subset L^2(H,\mu)$ is compact (see, e.g., [8]), we can construct a subsequence $\{\varphi_{n_k(\alpha)}^\alpha\}$ which is convergent in $L^2(H,\mu)$ and then another subsequence which is pointwise convergent. This implies that for each $\alpha \in (0,\frac{1}{2}]$, $\{\varphi_{n_k(\alpha)}\}$ is $\mu$-a.e. convergent on $K_\alpha = \{x : g(x) \le 1 - 2\alpha\}$.

Now, by a standard diagonal procedure we can find a subsequence $\{\varphi_{n_k}\}$ pointwisely convergent as required. $\quad\square$

2.5. *The Sobolev space* $W^{2,2}(K,\nu)$. It is easily seen that for all $h,k \in \mathbb{N}$ the linear operator

$$D_h D_k : C^2(K) \subset L^2(K,\nu) \to L^2(K,\nu), \qquad \varphi \mapsto D_h D_k \varphi,$$

is closable. If $\varphi$ belongs to the domain of the closure of $D_h D_k$ (which we shall still denote by $D_h D_k$) we shall say that $D_h D_k \varphi$ belongs to $L^2(K,\nu)$. Now we define $W^{2,2}(K,\nu)$ as the space of all functions $\varphi \in W^{1,2}(K,\nu)$ such that $D_h D_k \varphi \in L^2(K,\nu)$ for all $h,k \in \mathbb{N}$ and

$$\sum_{h,k=1}^\infty \int_H |D_h D_k \varphi(x)|^2 \nu(dx) < +\infty.$$

$W^{2,2}(K,\nu)$ is a Hilbert space with the inner product

$$\langle \varphi,\psi\rangle_{W^{2,2}(K,\nu)} = \langle \varphi,\psi\rangle_{W^{1,2}(K,\nu)} + \sum_{h,k=1}^\infty \int_K D_h D_k \varphi(x) D_h D_k \psi(x) \nu(dx).$$

If $\varphi \in W^{2,2}(K,\nu)$ we can define a Hilbert–Schmidt operator $D^2\varphi(x)$ on $K$ for $\nu$-almost all $x \in K$ by setting

$$\langle D^2\varphi(x)y,z\rangle = \sum_{h,k=1}^\infty D_h D_k \varphi(x)\langle y,e_h\rangle\langle z,e_k\rangle \qquad \forall y,z \in H.$$

We show now that if $\varphi \in W^{2,2}(K,\nu)$, then one can define the trace on $\Sigma$ of $D\varphi$. Similarly to the definition of the trace of $\varphi$ on $\Sigma$ we define $|Q^{1/2}\mathbf{n}(y)|\gamma(D\varphi) = \lim_{n\to\infty} |Q^{1/2}\mathbf{n}(y)|\gamma(D\varphi_N)$ in $L^2(\Sigma,\mu_\Sigma)$ for all $\{\varphi_n\} \subset C^2(K)$, $\varphi_n \to \varphi$ in $W^{2,2}(K,\nu)$.

Proposition 2.9 below shows that this trace is well defined.

PROPOSITION 2.9. *Assume that $\varphi \in W^{2,2}(K, \nu)$. Then:*

(i) $|Q^{1/2}\mathbf{n}(y)||\gamma(D\varphi)| \in L^2(\Sigma, \mu_\Sigma)$,

(ii) *the following estimate holds,*

$$(2.19) \qquad \int_\Sigma |Q^{1/2}\mathbf{n}(y)|^2 |\gamma(D\varphi(y))|^2 \mu_\Sigma(dy)$$
$$\leq C \left( \int_K |D\varphi(x)|^2 \nu(dx) + \int_K |\operatorname{Tr}[(D^2\varphi(x))^2]| \nu(dx) \right).$$

PROOF. Let $\varphi \in W^{2,2}(K, \nu)$ and let $\{\varphi_n\} \subset C^2(K)$ strongly convergent to $\varphi$ in $W^{2,2}(K, \nu)$. For $i \in \mathbb{N}$ we apply (2.15) to $D_i\varphi_n$. We have

$$\int_\Sigma |Dg(y)||Q^{1/2}\mathbf{n}(y)|^2 |D_i\varphi_n(y)|^2 \mu_\Sigma(dy)$$
$$\leq C \left( \int_K |D_i\varphi_n(x)|^2 \nu(dx) + \int_K |DD_i\varphi_n(x)|^2 \nu(dx) \right).$$

Summing up on $i$ yields

$$\int_\Sigma |Dg(y)||Q^{1/2}\mathbf{n}(y)|^2 |D\varphi_n(y)|^2 \mu_\Sigma(dy)$$
$$\leq C \left( \int_K |D\varphi_n(x)|^2 \nu(dx) + \sum_{i,j=1}^\infty \int_K |D_j D_i\varphi_n(x)|^2 \nu(dx) \right).$$

Then letting $n \to \infty$ we see that $\{Q^{1/2}\mathbf{n}(y)|\gamma(D\varphi_n)\}$ is strongly convergent in $L^2(K, \nu)$ and so (i) and (ii) follow. $\square$

When it will be no danger of confusion we shall simply set $D\varphi$ instead of $\gamma(D\varphi)$.

2.6. *The Sobolev space $W_A^{1,2}(K, \nu)$.* We define $W_A^{1,2}(K, \nu)$ as the space of all functions $\varphi \in W^{1,2}(K, \nu)$ such that

$$\sum_h^\infty \lambda_h \int_H |D_h\varphi(x)|^2 \nu(dx) < +\infty.$$

It is easy to see that $W_A^{1,2}(K, \nu)$ is a Hilbert space with the inner product

$$\langle \varphi, \psi \rangle_{W_A^{1,2}(K,\nu)} = \int_K \varphi(x)\psi(x)\nu(dx) + \sum_{h=1}^\infty \lambda_h \int_K D_h\varphi(x)D_h\psi(x)\nu(dx).$$

If $\varphi \in W_A^{1,2}(K, \nu)$ we can define an element of $K$, $A^{1/2}D\varphi(x)$ for $\nu$-almost all $x \in K$ by setting

$$\langle A^{1/2}D\varphi(x), y \rangle = \sum_{h=1}^\infty \lambda_h D_h\varphi(x)\langle y, e_h \rangle \qquad \forall y \in H.$$

**3. The Dirichlet form associated to $\nu$.** We define the symmetric Dirichlet form

$$a(\varphi, \psi) = \int_K \langle D\varphi, D\psi \rangle \, d\nu \qquad \forall \varphi, \psi \in D(a) = W^{1,2}(K, \nu) \times W^{1,2}(K, \nu).$$

Since, as seen earlier, $D$ is closed in $L^2(K, \nu)$ we infer that the form $a$ is closed in the sense of [15], page 315, and as a matter of fact the form $a$ is the closure of $a_0(\varphi, \psi) = \int_K \langle D\varphi, D\psi \rangle \, d\nu, \forall \varphi, \psi \in C_b^1(H)$.

By the Lax–Milgram theorem there exists an isomorphism

$$\mathcal{N} : W^{1,2}(K, \nu) \to (W^{1,2}(K, \nu))^*$$

[where $(W^{1,2}(K, \nu))^*$ is the dual space of $W^{1,2}(K, \nu)$] such that

$$\langle \varphi, \psi \rangle + a(\varphi, \psi) = \langle \mathcal{N}\varphi, \psi \rangle \qquad \forall \varphi, \psi \in W^{1,2}(K, \nu).$$

(Here $\langle \cdot, \cdot \rangle$ means the duality between $W^{1,2}(K, \nu)$ and $(W^{1,2}(K, \nu))^*$ which coincides with $\langle \cdot, \cdot \rangle_{L^2(K,\nu)}$ on $L^2(K, \nu)$.) We can identify $L^2(K, \nu)$ with its dual and, so, we have the well-known continuous and dense inclusions

$$W^{1,2}(K, \nu) \subset L^2(K, \nu) \subset (W^{1,2}(K, \nu))^*.$$

Now we define a linear operator $N : D(N) \subset L^2(K, \nu) \to L^2(K, \nu)$ as follows. We say that $\varphi \in D(N)$ if it belongs to $W^{1,2}(K, \nu)$ and that there exists $C > 0$ such that

$$(3.1) \qquad \left| \int_K \langle D\varphi, D\psi \rangle \, d\nu \right| \leq C |\psi|_{L^2(K,\nu)} \qquad \forall \psi \in W^{1,2}(K, \nu).$$

This inequality implies that $\mathcal{N}\varphi \in L^2(K, \nu)$. Finally, if $\varphi \in D(N)$ we set

$$N\varphi = \tfrac{1}{2}(I - \mathcal{N})\varphi.$$

In other words,

$$(3.2) \qquad \langle N\varphi, \psi \rangle = -\tfrac{1}{2} a(\varphi, \psi) \qquad \forall \varphi, \psi \in W^{1,2}(K, \nu).$$

THEOREM 3.1. *Operator $N$ is self adjoint in $L^2(K, \nu)$ and $\nu$ is an invariant measure for $N$,*

$$(3.3) \qquad \int_K N\varphi \, d\nu = 0 \qquad \forall \varphi \in D(N).$$

PROOF. By the closedness and symmetry of $a$ it follows that $N$ is closed and symmetric. Moreover, by the Lax–Milgram theorem, applied to symmetric bilinear form $(u, v) \to \lambda \langle u, v \rangle + a(u, v)$, we see that the range $R(\lambda I - N)$ of $\lambda I - N$ coincides with $L^2(K, \nu)$ for all $\lambda > 0$. Notice also that by (3.1)

$$(3.4) \qquad \langle N\varphi, \varphi \rangle = -\tfrac{1}{2} |D\varphi|^2_{L^2(K,\nu)} \qquad \forall \varphi \in D(N).$$

As regards (3.3) it is immediate by definition of $N$.   □

It is useful to notice also that for each $f \in L^2(K, \nu)$,

$$(\lambda I - N)^{-1} f = \{\varphi : \lambda \langle \varphi, \psi \rangle_{L^2(K,\nu)} + \tfrac{1}{2} a(\varphi, \psi)$$
$$= \langle f, \psi \rangle_{L^2(K,\nu)}, \forall \psi \in W^{1,2}(K, \nu)\}.$$

**4. The penalized problem.** We are here concerned with the penalized equation

$$(4.1) \qquad \begin{cases} dX_\varepsilon(t) + (AX_\varepsilon(t) + \beta_\varepsilon(X_\varepsilon(t))) \, dt = dW_t, \\ X_\varepsilon(0) = x, \end{cases}$$

where $\varepsilon > 0$, and

$$\beta_\varepsilon(x) = \frac{1}{\varepsilon}(x - \Pi_K(x)) \qquad \forall x \in H.$$

Since $\beta_\varepsilon$ is Lipschitz, (4.1) has a unique mild solution $X_\varepsilon(t, x)$.

The corresponding Kolmogorov operator reads as follows,

$$(4.2) \qquad N_\varepsilon \varphi = L\varphi - \langle \beta_\varepsilon(x), D\varphi \rangle, \qquad \varphi \in \mathcal{E}_A(H), \ \varepsilon > 0,$$

where $L$ is the Ornstein–Uhlenbeck operator

$$L\varphi = \tfrac{1}{2} \operatorname{Tr}[D^2 \varphi] - \langle x, AD\varphi \rangle, \qquad \varphi \in \mathcal{E}_A(H).$$

It is well known that $\nu_\varepsilon$ [defined in (2.1)–(2.3)] is an invariant measure for $N_\varepsilon$ and that

$$(4.3) \qquad \int_H N_\varepsilon \varphi \psi \, d\nu_\varepsilon = -\frac{1}{2} \int_H \langle D\varphi, D\psi \rangle \, d\nu_\varepsilon \qquad \forall \varphi, \psi \in \mathcal{E}_A(H).$$

Moreover, since $\beta_\varepsilon$ is Lipschitz continuous, operator $N_\varepsilon$ is essentially $m$-dissipative in $L^2(H, \nu_\varepsilon)$ (we still denote by $N_\varepsilon$ its closure) and $\mathcal{E}_A(H)$ is a core for $N_\varepsilon$ see [9].

Section 4.1 below is devoted to prove several estimates for the $(\lambda I - N_\varepsilon)^{-1} f$ where $f \in L^2(H, \nu_\varepsilon)$. Then these estimates are used in Section 4.2 to prove that $(\lambda I - N_\varepsilon)^{-1} f$ converges as $\varepsilon \to 0$ for any $f \in L^2(K, \nu)$ to $(\lambda I - N)^{-1} f$. Moreover we shall end up the section giving sharp informations about the domain of $N$.

4.1. *Estimates for $(\lambda I - N_\varepsilon)^{-1} f$.* We need a lemma.

LEMMA 4.1. *Let $\lambda > 0$, $\varphi \in \mathcal{E}_A(H)$ and set*

$$(4.4) \qquad f_\varepsilon = \lambda \varphi - N_\varepsilon \varphi.$$

*Then the following estimates hold*

$$(4.5) \qquad \int_H \varphi^2 \, d\nu_\varepsilon \le \frac{1}{\lambda^2} \int_H f_\varepsilon^2 \, d\nu_\varepsilon,$$

$$(4.6) \qquad \int_H |D\varphi|^2 \, d\nu_\varepsilon \le \frac{2}{\lambda} \int_H f_\varepsilon^2 \, d\nu_\varepsilon,$$

$$(4.7) \quad \begin{aligned} &\lambda \int_H |D\varphi|^2 \, d\nu_\varepsilon + \frac{1}{2} \int_H \mathrm{Tr}[(D^2\varphi)^2] \, d\nu_\varepsilon + \int_H |A^{1/2} D\varphi|^2 \, d\nu_\varepsilon \\ &\quad + \frac{1}{\varepsilon} \int_{K^c} \langle (I - D\Pi_K(x)) D\varphi, D\varphi \rangle \nu_\varepsilon \le 4 \int_H f_\varepsilon^2 \, d\nu_\varepsilon. \end{aligned}$$

PROOF. Multiplying both sides of (4.4) by $\varphi$, taking into account (4.3) and integrating in $\nu_\varepsilon$ over $H$, yields

$$(4.8) \qquad \lambda \int_H \varphi^2 \, d\nu_\varepsilon + \frac{1}{2} \int_H |D\varphi|^2 \, d\nu_\varepsilon = \int_H \varphi f_\varepsilon \, d\nu_\varepsilon.$$

Now (4.5) and (4.6) follow easily from the Hölder inequality. To prove (4.7) let us differentiate in the direction of $e_k$ both sides of (4.4). We obtain

$$\lambda D_k \varphi - N_\varepsilon D_k \varphi + \alpha_k D_k \varphi + \frac{1}{\varepsilon} \sum_{h=1}^\infty (\delta_{h,k} - \langle \Pi_K(x) e_h, e_k \rangle) D_h \varphi = D_k f_\varepsilon.$$

Multiplying both sides of (4.4) by $D_k \varphi$, taking into account (4.3), integrating in $\nu_\varepsilon$ over $H$ and then summing up over $k$, yields

$$(4.9) \quad \begin{aligned} &\lambda \int_H |D\varphi|^2 \, d\nu_\varepsilon + \frac{1}{2} \int_H \mathrm{Tr}[(D^2\varphi)^2] \, d\nu_\varepsilon + \int_H |A^{1/2} D\varphi|^2 \, d\nu_\varepsilon \\ &\quad + \frac{1}{\varepsilon} \int_{K^c} \langle (I - D\Pi_K(x)) D\varphi, D\varphi \rangle \, d\nu_\varepsilon = \int_H \langle D\varphi, Df_\varepsilon \rangle \, d\nu_\varepsilon. \end{aligned}$$

Noting finally that, again in view of (4.3),

$$\int_H \langle D\varphi, Df_\varepsilon \rangle \, d\nu_\varepsilon = 2 \int_H f_\varepsilon^2 \, d\nu_\varepsilon - 2\lambda \int_H f_\varepsilon \varphi \, d\nu_\varepsilon \le 4 \int_H f_\varepsilon^2 \, d\nu_\varepsilon,$$

the conclusion follows. □

Now we are able to prove the announced estimates.

PROPOSITION 4.2. *Let $\lambda > 0$, $f \in L^2(H, \nu_\varepsilon)$ and let $\varphi_\varepsilon$ be the solution of the equation*

$$(4.10) \qquad \lambda \varphi_\varepsilon - N_\varepsilon \varphi_\varepsilon = f.$$

*Then $\varphi_\varepsilon \in W^{2,2}(H, \nu_\varepsilon)$, $A^{1/2} D\varphi_\varepsilon \in L^2(H, \nu_\varepsilon)$ and the following estimates hold*

$$(4.11) \qquad \int_H \varphi_\varepsilon^2 \, d\nu_\varepsilon \le \frac{1}{\lambda^2} \int_H f^2 \, d\nu_\varepsilon,$$

$$(4.12) \qquad \int_H |D\varphi_\varepsilon|^2 \, d\nu_\varepsilon \le \frac{2}{\lambda} \int_H f^2 \, d\nu_\varepsilon,$$

$$(4.13) \qquad \begin{aligned} &\lambda \int_H |D\varphi_\varepsilon|^2 \, d\nu_\varepsilon + \frac{1}{2} \int_H \mathrm{Tr}[(D^2 \varphi_\varepsilon)^2] \, d\nu_\varepsilon + \int_H |A^{1/2} D\varphi_\varepsilon|^2 \, d\nu_\varepsilon \\ &\quad + \frac{1}{\varepsilon} \int_{K^c} \langle (I - D\Pi_K(x)) D\varphi_\varepsilon, D\varphi_\varepsilon \rangle \, d\nu_\varepsilon \le 4 \int_H f^2 \, d\nu_\varepsilon. \end{aligned}$$

PROOF. Inequality (4.11) is obvious since $N_\varepsilon$ is dissipative. Let us prove (4.12). Let $\lambda > 0$, $f \in L^2(H, \nu_\varepsilon)$ and let $\varphi_\varepsilon$ be the solution of (4.10). Since $\mathcal{E}_A(H)$ is a core for $N_\varepsilon$ there exists a sequence $\{\varphi_{\varepsilon,n}\}_{n \in \mathbb{N}} \subset \mathcal{E}_A(H)$ such that

$$\lim_{n \to \infty} \varphi_{\varepsilon,n} \to \varphi_\varepsilon, \qquad \lim_{n \to \infty} N_\varepsilon \varphi_{\varepsilon,n} \to N_\varepsilon \varphi_\varepsilon \qquad \text{in } L^2(H, \nu_\varepsilon).$$

Set $f_{\varepsilon,n} = \lambda \varphi_{\varepsilon,n} - N_\varepsilon \varphi_{\varepsilon,n}$. Clearly, $f_{\varepsilon,n} \to f$ as $n \to \infty$ in $L^2(H, \nu_\varepsilon)$.

We claim that $\varphi_\varepsilon \in W^{1,2}(H, \nu_\varepsilon)$ and that

$$\lim_{n \to \infty} D\varphi_{\varepsilon,n} \to D\varphi_\varepsilon \qquad \text{in } L^2(H, \nu_\varepsilon; H).$$

Let in fact $m, n \in \mathbb{N}$, then by (4.6) it follows that

$$\int_H |D\varphi_{\varepsilon,n} - D\varphi_{\varepsilon,m}|^2 \, d\nu_\varepsilon \le \frac{1}{\lambda^2} \int_H |f_{\varepsilon,n} - f_{\varepsilon,m}|^2 \, d\nu_\varepsilon.$$

Therefore the sequence $\{\varphi_{\varepsilon,n}\}_{n \in \mathbb{N}}$ is Cauchy in $W^{1,2}(H, \nu_\varepsilon)$ and the claim follows. Estimate (4.13) can be proved similarly. $\square$

We conclude this section with an integration by parts formula needed later. We set

$$V = \{\psi \in W^{1,2}(K, \nu) : |Q^{1/2} \mathbf{n}| \psi \in L^2(\Sigma, \mu_\Sigma)\}.$$

LEMMA 4.3. *Let $\varphi \in D(N_\varepsilon)$ and $\psi \in V$. Then the following identity holds.*

$$(4.14) \qquad \begin{aligned} \int_K N_\varepsilon \varphi \psi \, d\nu = &-\frac{1}{2} \int_K \langle D\varphi, D\psi \rangle \, d\nu \\ &+ \frac{1}{\mu(K)} \int_\Sigma \langle \gamma(D\varphi), \mathbf{n}(y) \rangle \psi \, d\mu_\Sigma. \end{aligned}$$

PROOF.  Taking in account that $\mathcal{E}_A(H)$ is a core for $N_\varepsilon$, it is sufficient to prove (4.14) for $\varphi \in \mathcal{E}_A(H)$. By the basic integration by parts formula we deduce, for any $i \in \mathbb{N}$ and $\psi \in V$ that

$$\int_K D_i\varphi D_i\psi \, d\nu = -\int_K D_i^2\varphi\psi \, d\nu + \frac{1}{\mu(K)}\int_\Sigma \gamma(D_i\varphi)(\mathbf{n}(y))_i\psi \, d\mu_\Sigma$$

$$+ \frac{1}{\lambda_i}\int_K x_i D_i\varphi\psi \, d\nu.$$

Now, summing up on $i$ yields

$$\int_K \langle D\varphi, D\psi\rangle \, d\nu = -\int_K \mathrm{Tr}[D^2\varphi]\psi \, d\nu + \frac{1}{\mu(K)}\int_\Sigma\langle\gamma(D\varphi), \mathbf{n}(y)\rangle\psi \, d\mu_\Sigma$$

$$+ 2\int_K \langle x, AD\varphi\rangle\psi \, d\nu.$$

That is nothing else but (4.14).  □

4.2.  *Convergence of* $\{\varphi_\varepsilon\}$.  We are going to show that the sequence $\{\varphi_\varepsilon\}$ is convergent in $L^2(K, \nu)$. We first note that for $f \in C_b(H)$ we have

(4.15)          $$\varphi_\varepsilon(x) = \mathbb{E}\int_0^\infty e^{-\lambda t}f(X_\varepsilon(t, x)) \, dt \qquad \forall x \in H.$$

Now, by a standard argument it follows that from (4.15) that if $f \in C_b^1(H)$ we have

(4.16)          $$\sup_{x\in H}|D\varphi_\varepsilon(x)| \leq \frac{1}{\lambda}\|Df\|_{C_b(H)} \qquad \forall\varepsilon, \lambda > 0.$$

Theorem 4.4 is the main result of this section.

THEOREM 4.4.  *Let* $\lambda > 0$, $f \in L^2(K, \nu)$ *and let* $\varphi_\varepsilon$ *be the solution of* (4.10). *Then* $\{\varphi_\varepsilon\}$ *is strongly convergent in* $L^2(K, \nu)$ *to* $\varphi = (\lambda I - N)^{-1}f$ *where* $N$ *is defined by* (3.1).

*Moreover, the following statements hold:*

(i) $\lim_{\varepsilon\to 0} D\varphi_\varepsilon = D\varphi$ *in* $L^2(K, \nu; H)$,

(ii) $\varphi \in W_A^{1,2}(H, \nu) \cap W^{2,2}(K, \nu)$,

(iii) $\varphi$ *fulfills the Neumann condition*

(4.17)          $$\frac{d\varphi}{dn}(x) = \langle D\varphi(x), \mathbf{n}(x)\rangle = 0 \qquad on \ \Sigma,$$

*where* $\langle D\varphi(x), \mathbf{n}(x)\rangle$ *is defined by Proposition 2.9 and* $|Q^{1/2}\mathbf{n}(x)|\langle D\varphi(x), \mathbf{n}(x)\rangle \in L^2(\Sigma, \mu_\Sigma)$.

PROOF. Without danger of confusion we shall denote again by $f$ the restriction $f|_K$ of $f$ to $K$. In fact each $f \in L^2(K, \nu)$ can be extended by 0 outside $K$ to a function in $L^2(H, \nu)$. By this convention, everywhere in the sequel $(\lambda I - N)^{-1} f$ for $f \in L^2(H, \nu)$ means $(\lambda I - N)^{-1} f|_K$.

*Step* 1. We have

$$(4.18) \qquad \lim_{\varepsilon \to 0} \varphi_\varepsilon = (\lambda I - N)^{-1} f \qquad \text{in } L^2(K, \nu).$$

In fact by (4.11), (4.12) and the compactness of the embedding of $W^{1,2}(K, \nu)$ in $L^2(K, \nu)$ it follows that there exist a sequence $\{\varepsilon_k\} \to 0$ and $\varphi \in W^{1,2}(K, \nu)$ such that

$$\varphi_{\varepsilon_k} \to \varphi \qquad \text{strongly in } L^2(K, \nu),$$

$$D\varphi_{\varepsilon_k} \to D\varphi \qquad \text{weakly in } L^2(K, \nu).$$

Let $\psi \in C_b^1(H)$ and consider the identity

$$\frac{1}{2} \int_H \langle D\varphi_\varepsilon, D\psi \rangle \, d\nu_\varepsilon = \int_H (f - \lambda \varphi_\varepsilon) \psi \, d\nu_\varepsilon,$$

which is equivalent to

$$(4.19) \qquad \frac{1}{2} \int_K \langle D\varphi_\varepsilon, D\psi \rangle \, d\nu + \frac{1}{2} \int_{K^c} \langle D\varphi_\varepsilon, D\psi \rangle \, d\nu_\varepsilon = \int_H (f - \lambda \varphi_\varepsilon) \psi \, d\nu_\varepsilon.$$

Since, we have

$$\left| \int_{K^c} \langle D\varphi_\varepsilon, D\psi \rangle \, d\nu_\varepsilon \right|^2 \le \int_H |D\varphi_\varepsilon|^2 \, d\nu_\varepsilon \int_{K^c} |D\psi|^2 \, d\nu_\varepsilon$$

$$\le \frac{2}{\lambda} \int_H f^2 \, d\nu_\varepsilon \int_{K^c} |D\psi|^2 \, d\nu_\varepsilon \to 0$$

as $\varepsilon \to 0$, we deduce, letting $\varepsilon \to 0$ in (4.19) that

$$\frac{1}{2} \int_K \langle D\varphi, D\psi \rangle \, d\nu = \int_K (f - \lambda \varphi) \psi \, d\nu \qquad \forall \psi \in C_b^1(H).$$

Obviously, this identity extends to all $\psi \in W^{1,2}(H, \nu)$, which implies that $\varphi = (\lambda I - N)^{-1} f$ and that $\varphi_\varepsilon \to \varphi$ strongly in $L^2(K, \nu)$.

*Step* 2. We have

$$\lim_{\varepsilon \to 0} D\varphi_\varepsilon = D\varphi \qquad \text{in } L^2(K, \nu; K).$$

We first assume that $f \in C_b^1(H)$. Let us start from the identity (4.8),

$$(4.20) \qquad \frac{1}{2} \int_H |D\varphi_\varepsilon|^2 \, d\nu_\varepsilon = \int_K (\lambda \varphi_\varepsilon - f) \varphi_\varepsilon \, d\nu_\varepsilon,$$

which implies

$$(4.21) \quad \lim_{\varepsilon \to 0} \frac{1}{2} \int_H |D\varphi_\varepsilon|^2 \, d\nu_\varepsilon = \int_K (\lambda\varphi - f)\varphi \, d\nu$$

$$= -\langle N\varphi, \varphi \rangle = \frac{1}{2} \int_K |D\varphi|^2 \, d\nu.$$

Here we have used the fact that

$$\lim_{\varepsilon \to 0} \int_{K^c} |D\varphi_\varepsilon|^2 \, d\nu_\varepsilon(x) = 0,$$

which follows taking into account (4.16).

Therefore there exists a sequence $\{\varepsilon_k\}$ such that

$$\varphi_{\varepsilon_k} \to \varphi, \qquad \text{strongly in } L^2(K, \nu),$$

$$D\varphi_{\varepsilon_k} \to D\varphi, \qquad \text{weakly in } L^2(K, \nu; H),$$

$$\lim_{k \to \infty} \int_K |D\varphi_{\varepsilon_k}|^2 \, d\nu = \int_K |D\varphi|^2 \, d\nu.$$

This implies that $D\varphi_{\varepsilon_k} \to D\varphi$ strongly in $L^2(K, \nu; H)$.

We finally assume that $f \in L^2(H, \nu)$. Since $C_b^1(H)$ is dense in $L^2(K, \nu)$, there exists a sequence $\{f_n\} \subset C_b^1(H)$ strongly convergent in $L^2(K; \nu)$ to $f$. Set $\varphi_{n,\varepsilon} = (\lambda I - N_\varepsilon)^{-1} f_n$. By (4.12) we have

$$\int_H |D\varphi_\varepsilon - D\varphi_{n,\varepsilon}|^2 \, d\nu_\varepsilon \le \frac{2}{\lambda} \int_K |f - f_n|^2 \, d\nu,$$

which implies

$$\int_K |D\varphi_\varepsilon - D\varphi_{n,\varepsilon}|^2 \, d\nu \le \frac{2}{\lambda} \int_K |f - f_n|^2 \, d\nu.$$

So, again $D\varphi_{\varepsilon_k} \to D\varphi$ strongly in $L^2(K, \nu; H)$.

*Step* 3. We have

$$(4.22) \quad \varphi \in W_A^{1,2}(K, \nu; H) \cap W^{2,2}(K; \nu).$$

By estimate (4.13) we have that $\{\varphi_\varepsilon\}$ is bounded in $W^{2,2}(K, \nu)$. Therefore there is a subsequence, still denoted $\{\varphi_\varepsilon\}$ which converges to $\varphi$ in $W^{2,2}(K, \nu)$. In the same way we see that $\varphi \in W_A^{1,2}(K, \nu; H)$.

*Step* 4. Checking the Neumann condition for $\varphi$.

From (4.14) we get

$$\int_K N_\varepsilon \varphi_\varepsilon \psi \, d\nu = -\frac{1}{2} \int_K \langle D\varphi_\varepsilon, D\psi \rangle \, d\nu + \frac{1}{\mu(K)} \int_\Sigma \psi \langle \gamma(D\varphi_\varepsilon), \mathbf{n}(y) \rangle \, d\mu_\Sigma.$$

Recalling that $N_\varepsilon \varphi_\varepsilon = \lambda \varphi_\varepsilon - f \longrightarrow \lambda \varphi - f = N\varphi$ in $L^2(K, \nu)$ and that $|Q^{1/2}\mathbf{n}(y)|\langle \gamma(D\varphi_\varepsilon), \mathbf{n}(y)\rangle \to |Q^{1/2}\mathbf{n}(y)|\langle \gamma(D\varphi), \mathbf{n}(y)\rangle$ in $L^2(\Sigma, \mu_\Sigma)$ by Proposition 2.9, by (i) and by (3.4) we obtain

$$\int_\Sigma \langle \gamma(D\varphi), \mathbf{n}(y)\rangle \psi \, d\mu_\Sigma = 0 \qquad \forall \psi \in V,$$

which implies (4.17) as claimed. [The set $\{\gamma(\psi) : \psi \in V\}$ is dense in $L^2(\Sigma, \mu_\Sigma)$.] This completes the proof. $\quad\square$

In particular, taking into account that $D(N)$ is equal to the range of $(\lambda I - N)^{-1}$ we derive by Theorem 4.4 the following result, which gives a sharp information on the structure of the domain of $N$.

COROLLARY 4.5. *We have*

$$(4.23) \quad D(N) \subset \left\{ \varphi \in W^{1,2}_A(H, \nu) \cap W^{2,2}(K, \nu) : \frac{d}{dn}\varphi(x) = 0 \ on \ \Sigma \right\}.$$

We notice also that for $\varphi \in D(N)$ regular $N\varphi$ is the classical elliptic differential operator in $H$. More precisely, we have

COROLLARY 4.6. *If* $\operatorname{Tr} D^2\varphi \in L^2(K, \nu)$, $\langle x, AD\varphi \rangle \in L^2(K, \nu)$ *and* $\frac{d\varphi}{dn}(y) = 0, \forall y \in \Sigma$ *then* $\varphi \in D(N)$ *and*

$$(4.24) \qquad N\varphi(x) = \tfrac{1}{2}\operatorname{Tr} D^2\varphi - \langle x, AD\varphi \rangle \qquad \forall x \in \overset{\circ}{K}.$$

PROOF. By integration by parts formula (2.11) we see that

$$(4.25) \quad \begin{aligned} \int_K \langle D\varphi, D\psi \rangle \nu(dx) &= -\int_K \left( \tfrac{1}{2}\operatorname{Tr} D^2\varphi - \langle x, AD\varphi \rangle \right) \nu(dx) \\ &\quad + \frac{1}{\mu(K)}\int_\Sigma \psi(y)\frac{d\varphi}{dn}(y)\mu_\Sigma(dy) \qquad \forall \psi \in V, \end{aligned}$$

which in virtue of (iv) and (3.2) implies (4.24) as claimed. $\quad\square$

REMARK 4.7. We conjecture that in Corollary 4.5 one has equality in relation (4.23), but we failed to prove it. This happens when $N$ is replaced by the Ornstein–Uhlenbeck generator $L$ and $\nu$ by the Gaussian measure $\mu$ (see [9]).

Notice also that if $\varphi \in D(N)$ we cannot conclude that $\operatorname{Tr} D^2\varphi \in L^2(K, \nu)$ and $\langle x, AD\varphi \rangle \in L^2(K, \nu)$. This is obviously true if $H$ is finite-dimensional.

## 5. Perturbation results.

5.1. *Perturbation by a regular gradient.* Let us consider the stochastic differential inclusion,

$$(5.1) \quad \begin{cases} dX(t) + (AX(t) + DV(X(t)) + N_K(X(t))) \, dt \ni dW(t), \\ X(0) = x, \end{cases}$$

where $A, K$ and $W$ are as before and $V : H \to \mathbb{R}$ is a $C^2$ function such that $DV \in C_b^1(H; H)$.

Let us introduce a probability measure $\zeta \in \mathcal{P}(K)$ by setting

$$\zeta(dx) = Z_\zeta^{-1} e^{-2V(x)} \nu(dx),$$

where

$$Z_\zeta = \int_K e^{-2V(y)} \nu(dy).$$

Arguing as in the proof of (2.11), we can show the following integration by parts formula.

THEOREM 5.1. *Let $\varphi \in C_b^1(H)$. Then for any $z \in H$ we have*

$$
(5.2) \quad
\begin{aligned}
\int_K &\langle D\varphi(x), Q^{1/2} z\rangle \zeta(dx) \\
&= \int_K \varphi(x)\langle DV(x), Q^{1/2} z\rangle \zeta(dx) \\
&\quad + \frac{1}{2\mu(K)Z_\zeta} \int_\Sigma \varphi(y)\langle \mathbf{n}(y), Q^{1/2} z\rangle e^{-2U(y)} \mu_\Sigma(dy) \\
&\quad + \int_K W_z(x)\varphi(x)\zeta(dx).
\end{aligned}
$$

Now all considerations of Sections 2, 3 and 4 can be easily generalized. In particular, estimate (4.7) reads as follows

$$
(5.3) \quad
\begin{aligned}
\lambda &\int_H |D\varphi|^2 \, d\zeta_\varepsilon + \frac{1}{2} \int_H \mathrm{Tr}[(D^2\varphi)^2] \, d\zeta_\varepsilon + \int_H |A^{1/2} D\varphi|^2 \, d\zeta_\varepsilon \\
&+ \int_H \langle D^2 V \cdot D\varphi, D\varphi\rangle \, d\zeta_\varepsilon + \frac{1}{\varepsilon} \int_{K^c} \langle (I - D\Pi_K(x))D\varphi, D\varphi\rangle \, d\zeta_\varepsilon \\
&\leq 4 \int_H f_\varepsilon^2 \, d\zeta_\varepsilon.
\end{aligned}
$$

In conclusion, we arrive at the following result.

THEOREM 5.2. *The operator $N$ (defined as in Section 3 with the Dirichlet form induced by $\zeta$) is self adjoint in $L^2(K, \zeta)$ and $\zeta$ is an invariant measure for $N$,*

$$(5.4) \qquad \int_K N\varphi \, d\zeta = 0 \qquad \forall \varphi \in D(N).$$

*Moreover, we have*

$$(5.5) \qquad D(N) \subset \left\{ \varphi \in W_A^{1,2}(H, \zeta) \cap W^{2,2}(K, \zeta) \colon \frac{d}{dn} \varphi(x) = 0 \ on \ \Sigma \right\}.$$

(Details are omitted.)

5.2. *Perturbation by a bounded Borel drift.* Let $F \colon H \to H$ be bounded and Borel and consider the stochastic differential inclusion,

$$(5.6) \qquad \begin{cases} dX(t) + (AX(t) + F(X(t)) + N_K(X(t))) \, dt \ni dW(t), \\ X(0) = x. \end{cases}$$

Let moreover $G$ be the linear operator in $L^2(K, \nu)$ defined as

$$(5.7) \qquad G\varphi = N\varphi + \langle F(x), D\varphi \rangle, \qquad \varphi \in D(N).$$

PROPOSITION 5.3. *$G$ is the infinitesimal generator of a strongly continuous compact semigroup $Q_t$ on $L^2(K, \nu)$. Moreover its resolvent $(\lambda I - G)^{-1}$ is given by*

$$(5.8) \qquad (\lambda I - G)^{-1} = (\lambda I - N)^{-1}(1 - T_\lambda)^{-1}, \qquad \lambda > \lambda_0,$$

*where*

$$(5.9) \qquad \lambda_0 = 2\|F\|_0^2 = 2 \sup_{x \in H} |F(x)|^2$$

*and*

$$(5.10) \qquad T_\lambda \psi(x) = \langle F(x), D(\lambda I - N)^{-1} \psi(x) \rangle, \qquad \psi \in L^2(K, \nu), \ x \in K.$$

PROOF. Let $\lambda > 0$, $f \in L^2(K, \nu)$. Consider the equation

$$(5.11) \qquad \lambda\varphi - N\varphi - \langle F(x), D\varphi \rangle = f.$$

Setting $\psi = \lambda\varphi - N\varphi$ (5.11) becomes

$$(5.12) \qquad \psi - T_\lambda \psi = f,$$

where $T_\lambda$ is defined by (5.10).

On the other hand, by (4.12) it follows that

$$\int_H |D(\lambda I - N)^{-1} \psi|^2 \, d\nu \le \frac{2}{\lambda} \int_H \psi^2 \, d\nu,$$

so that

$$\|T_\lambda \psi\|_{L^2(H,\mu)} \leq \sqrt{\frac{2}{\lambda}} \|F\|_0 \|\psi\|_{L^2(H,\mu)}.$$

Therefore if $\lambda > \lambda_0$ (5.11) has a unique solution and the conclusion follows.

Finally, the compactness property of $Q_t$ for $t > 0$ follows from (5.9) and the compactness of operator $(\lambda I - N)^{-1}$. $\square$

We want now to show that operator $G$ possesses an invariant measure $\zeta$ absolutely continuous with respect to $\nu$. For this let us consider the adjoint semigroup $Q_t^*$; we denote by $G^*$ its infinitesimal generator, and by $\Sigma^*$ the set of all its stationary points:

$$\Sigma^* = \{\varphi \in L^2(K,\nu) : Q_t^* \varphi = \varphi, t \geq 0\}.$$

Though the following lemma is standard, we give a proof, however, for reader's convenience. We shall denote by $\mathbb{1}$ the functions identically equal to 1.

LEMMA 5.4. *$Q_t^*$ has the following properties:*

    (i) *For all $\varphi \geq 0$ $\nu$-a.e., one has $Q_t^* \varphi \geq 0$ $\nu$-a.e.*

    (ii) *$\Sigma^*$ is a lattice, that is, if $\varphi \in \Sigma^*$ then $|\varphi| \in \Sigma^*$.*

PROOF. Let $\psi_0 \geq 0$ $\nu$-a.e. Then for all $\varphi \geq 0$ $\nu$-a.e. and all $t > 0$ we have

$$\int_K Q_t \varphi \psi_0 \, d\nu = \int_K \varphi Q_t^* \psi_0 \, d\nu \geq 0.$$

This implies that $\psi_0 \geq 0$ $\nu$-a.e., and (i) is proved.

Let us prove (ii). Assume that $\varphi \in \Sigma^*$, so that $\varphi(x) = Q_t^* \varphi(x)$. Then we have

$$(5.13) \qquad |\varphi(x)| = |Q_t^* \varphi(x)| \leq Q_t^*(|\varphi|)(x).$$

We claim that

$$|\varphi(x)| = Q_t^*(|\varphi|)(x), \qquad x - \nu \text{ a.s.}$$

Assume by contradiction that there is a Borel subset $I \subset K$ such that $\nu(I) > 0$ and

$$|\varphi(x)| < Q_t^*(|\varphi|)(x) \qquad \forall x \in I.$$

Then we have

$$(5.14) \qquad \int_K |\varphi(x)| \nu(dx) < \int_K Q_t^*(|\varphi|)(x)\nu(dx).$$

On the other hand,

$$\int_K Q_t^*(|\varphi|)\,d\mu = \langle Q_t^*(|\varphi|), \mathbb{1}\rangle_{L^2(K,\nu)} = \langle|\varphi|, \mathbb{1}\rangle_{L^2(K,\mu)} = \int_K |\varphi|\,d\mu,$$

which is in contradiction with (5.14).  □

The following result is a generalization of a similar result concerning the Ornstein–Uhlenbeck semigroup proved in [8].

PROPOSITION 5.5.  *There exists an invariant measure $\zeta$ of $Q_t$ which is absolutely continuous with respect to $\nu$. Moreover*

$$\rho := \frac{d\zeta}{d\nu} \in L^2(K,\nu).$$

PROOF.  Let $\lambda > 0$ be fixed. Clearly $\mathbb{1} \in D(G)$ and we have $G\mathbb{1} = 0$. Consequently $\frac{1}{\lambda}$ is an eigenvalue of $(\lambda I - G)^{-1}$ since

$$(\lambda I - G)^{-1}\mathbb{1} = \frac{1}{\lambda}\mathbb{1}.$$

Moreover $\frac{1}{\lambda}$ has finite multiplicity because $(\lambda I - G)^{-1}$ is compact. Therefore $((\lambda I - G)^{-1})^*$ is compact as well and $\frac{1}{\lambda}$ is an eigenvalue for $((\lambda I - G)^{-1})^*$. Consequently there exists $\rho \in L^2(K,\nu)$ not identically equal to zero such that

$$(5.15)\qquad ((( \lambda I - G)^{-1})^*)\rho = \frac{1}{\lambda}\rho.$$

It follows that $\rho \in D(G)$ and $G^*\rho = 0$. Since $\Sigma^*$ is a lattice, $\rho$ can be chosen to be nonnegative and such that $\int_K \rho\,d\nu = 1$.

Now set

$$\zeta(dx) = \rho(x)\nu(dx), \qquad x \in K.$$

We claim that $\zeta$ is an invariant measure for $Q_t$. In fact taking the inverse Laplace transform in (5.15) we find that

$$Q_t^*\rho = \rho,$$

which implies that for any $\varphi \in L^2(K,\nu)$,

$$\int_K Q_t\varphi\,d\zeta = \int_K Q_t\varphi\rho\,d\nu = \int_K \varphi Q_t^*\rho\,d\nu = \int_K \varphi\,d\zeta.$$

The proof is complete.  □

Notice now that, since $\frac{d\zeta}{d\nu} \in L^2(K,\nu)$ there is a natural inclusion of $L^2(K,\nu)$ in $L^1(K,\zeta)$ so, we can introduce the linear operator in $L^1(K,\zeta)$,

$$(5.16)\qquad N_F : D(N) \subset L^2(K,\nu) \to L^1(K,\zeta), \qquad N_F\varphi := G\varphi.$$

This is the final result of the paper.

PROPOSITION 5.6. *Operator $N_F$ defined by (5.16) is dissipative in $L^1(K, \zeta)$ and its closure is $m$-dissipative.*

PROOF. The dissipativity of operator $N_F$ in $L^1(K, \zeta)$ follows from the fact that measure $\zeta$ is invariant for $N_F$ and a standard argument; see [11]. Moreover the range of $\lambda I - N_F$ contains $L^2(K, \nu)$ for $\lambda > \lambda_0$ which is dense in $L^1(K, \zeta)$. So, the conclusion follows from the Lumer–Phillips theorem. $\square$

## 6. An example.

EXAMPLE 6.1. Consider the stochastic equation

$$(6.1) \quad \begin{cases} dX(t) - \Delta X(t)\, dt + N_K(X(t))\, dt \ni dW_t, & \text{in } (0, \infty) \times \mathcal{O}, \\ X(t) = 0, & \text{on } \partial \mathcal{O}, \\ X(0) = x, & \text{in } \mathcal{O}, \end{cases}$$

where $\mathcal{O}$ is a bounded and open interval of $\mathbb{R}$, and

$$K = \{x \in L^2(\mathcal{O}) : \|x\|_{L^2(\mathcal{O})} \leq \rho\}.$$

Then the previous results apply with $H = L^2(\mathcal{O})$, $A = -\Delta$, $D(A) = H_0^1(\mathcal{O}) \cap H^2(\mathcal{O})$.

Thus the Markov semigroup $P_t$ generated by $N$ in this case is given by $(P_t \varphi_0)(x) = \varphi(t, x)$ where

$$\varphi \in C^1([0, \infty); L^2(L^2(\mathcal{O}), \nu)) \cap C([0, \infty); W^{2,2}(K, \nu)) \cap W_A^{1,2}(K, \nu; L^2(\mathcal{O}))$$

is the solution to infinite-dimensional parabolic equation

$$\begin{cases} \dfrac{d}{dt} \displaystyle\int_K \varphi(t, x) \psi(x) \nu(dx) \\ \qquad + \displaystyle\int_K \left( \int_{\mathcal{O}} D\varphi(t, x)(\xi) D\psi(x)(\xi)\, d\xi \right) \nu(dx), & \forall t \geq 0, \\ \varphi(0, x) = \varphi_0(x), & x \in L^2(\mathcal{O}) \end{cases}$$

for all $\psi \in W^{1,2}(K, \nu)$.

A more general case is that where

$$(6.2) \quad K = \left\{ x \in L^2(\mathcal{O}) : \int_{\mathcal{O}} j(x(\xi))\, d\xi \leq \rho^2 \right\},$$

where $j : \mathbb{R} \to \mathbb{R}$ is a $C^\infty$ function such that $0 \leq j(r) \leq C_1 r^2$, $j''(r) \geq C_2 > 0$, $\forall r \in \mathbb{R}$. In this latter case

$$\Sigma = \left\{ x : \int_{\mathcal{O}} j(x(\xi))\, d\xi = \rho^2 \right\} \quad \text{and} \quad N_K(x)(\xi) = \{\lambda \nabla j(x(\xi))\}_{\lambda > 0} \qquad \forall x \in \Sigma.$$

## APPENDIX

Here we shall present for the reader's convenience a few results on co-aerea formula used in Section 2.1, under additional conditions on $g$, Hypothesis A.1.

**A.1. The co-area formula.** Let $H$ be a separable Hilbert space and $\mu = N_Q$ a nondegenerate Gaussian measure in $H$. Let $(e_k)$ be the complete orthonormal basis in $H$ corresponding to the eigenvalues $(\lambda_k)$, a sequence of positive numbers, that is, $Qe_k = \lambda_k e_k, k \in \mathbb{N}$.

Let us recall the integration by parts formula

$$(A.1) \qquad \int_H D_h \varphi \psi \, d\mu = -\int_H D_h \psi \varphi \, d\mu + \frac{1}{\lambda_h} \int_H x_h \varphi \psi \, d\mu$$

for any $\varphi$ bounded and Borel.

We are given a Borel bounded mapping $g : H \to \mathbb{R}$ of class $C^2$ such that

HYPOTHESIS A.1.

$$
I_1 := \int_H \frac{\mathrm{Tr}[Q D^2 g(x)]}{|Q^{1/2} Dg(x)|^2} \mu(dx) < \infty,
$$

$$
(A.2) \qquad I_2 := \int_H \frac{\langle D^2 g(x) \cdot Q^{1/2} g(x), Q^{1/2} g(x)\rangle}{|Q^{1/2} Dg(x)|^4} \mu(dx) < \infty,
$$

$$
I_3 := \int_H \frac{\langle x, Dg(x)\rangle}{|Q^{1/2} Dg(x)|^2} \mu(dx) < \infty.
$$

REMARK A.2. Let $\rho$ be a nonnegative $C^2$ real function such that for some $c > 0$, $m \in \mathbb{N}$

$$|\rho'(r)| + |\rho''(r)| \le c(1 + r^m), \qquad |\rho'(r)| \ge c$$

and set $g(x) = \rho(|x|^2)$. Then we have

$$
I_1 = \frac{2\rho''(|x|^2)|Q^{1/2}x|^2 + \rho'(|x|^2)\,\mathrm{Tr}\,Q}{\rho'(|x|^2)|Q^{1/2}x|^2},
$$

$$
I_2 = \frac{4(2\rho''(|x|^2) + \rho'(|x|^2)}{\rho'(|x|^2)|Q^{1/2}x|^2},
$$

$$
I_3 = \frac{|x|^2}{\rho'(|x|^2)|Q^{1/2}x|^2}.
$$

Then it is not difficult to see that Hypothesis A.1 is fulfilled. Let us check for instance that

$$(A.3) \qquad J_1 : \int_H \frac{1}{|Q^{1/2}x|^2} \mu(dx) < +\infty.$$

We have in fact

$$\frac{1}{|Q^{1/2}x|^2} = \int_0^{+\infty} e^{-t|Q^{1/2}x|^2}\, dt,$$

so that

$$J_1 = \int_0^{+\infty} dt \int_H e^{-t|Q^{1/2}x|^2} \mu(dx) = \int_0^{+\infty} \prod_{k=1}^{\infty} \frac{1}{\sqrt{1 + 2t\lambda_k^2}}\, dt$$

$$\leq \int_0^{+\infty} \prod_{k=1}^{3} \frac{1}{\sqrt{1 + 2t\lambda_k^2}}\, dt < +\infty.$$

We denote by $\mu_g := g_\# \mu$ the law of $g$ on $(\mathbb{R}, \mathcal{B}(\mathbb{R}))$. Then for any $\varphi \colon \mathbb{R} \to \mathbb{R}$ it holds

$$(A.4) \qquad \int_{\mathbb{R}} \varphi(r) \mu_g(dr) = \int_H \varphi(g(x)) \mu(dx).$$

We are going to show following [4] that, under Hypothesis A.1, $g_\# \mu \ll \ell$, where $\ell$ is the Lebesgue measure on $\mathbb{R}$, using the well-known sufficient condition

$$(A.5) \qquad \left| \int_{\mathbb{R}} \varphi'(r) \mu_g(dr) \right| = \left| \int_H \varphi'(g(x)) \mu(dx) \right| \leq C \|\varphi\|_0 \qquad \forall \varphi \in C_b^1(H).$$

PROPOSITION A.3. *Assume that Hypothesis A.1 is fulfilled. Then* $g_\# \mu = \mu_g \ll \ell$.

PROOF. We claim that

$$(A.6)$$
$$\int_H \varphi'(g(x)) \mu(dx)$$
$$= -\int_H \varphi(g(x)) \frac{\operatorname{Tr}[Q D^2 g(x)]}{|Q^{1/2} Dg(x)|^2} \mu(dx)$$
$$\quad - 2 \int_H \varphi(g(x)) \frac{\langle D^2 g(x) \cdot Q^{1/2} g(x), Q^{1/2} g(x) \rangle}{|Q^{1/2} Dg(x)|^4} \mu(dx)$$
$$\quad + \int_H \varphi(g(x)) \frac{\langle x, Dg(x) \rangle}{|Q^{1/2} Dg(x)|^2} \mu(dx),$$

which will yield the conclusion.

Since

$$\langle D\varphi(g(x)), Q Dg(x) \rangle = \varphi'(g(x)) |Q^{1/2} Dg(x)|^2,$$

we have

$$\int_H \varphi'(g(x))\mu(dx) = \int_H \frac{1}{|Q^{1/2}Dg(x)|^2}\langle D\varphi(g(x)), QDg(x)\rangle\mu(dx)$$

$$= \sum_{k=1}^\infty \lambda_k \int_H D_k\varphi(g(x))\frac{D_k g(x)}{|Q^{1/2}Dg(x)|^2}\mu(dx).$$

Using (A.1) yields

$$\int_H \varphi'(g(x))\mu(dx) = \sum_{k=1}^\infty \lambda_k \int_H D_k\varphi(g(x))\frac{D_k g(x)}{|Q^{1/2}Dg(x)|^2}\mu(dx)$$

$$= -\sum_{k=1}^\infty \lambda_k \int_H \varphi(g(x))D_k\left[\frac{D_k g(x)}{|Q^{1/2}Dg(x)|^2}\right]\mu(dx)$$

$$+ \sum_{k=1}^\infty \int_H x_k\varphi(g(x))\frac{D_k g(x)}{|Q^{1/2}Dg(x)|^2}\mu(dx).$$

But

$$D_k\left[\frac{D_k g(x)}{|Q^{1/2}Dg(x)|^2}\right] = \frac{D_k^2 g(x)}{|Q^{1/2}Dg(x)|^2} - 2\frac{\sum_{j=1}^\infty \lambda_j D_k D_j g(x)D_j g(x)}{|Q^{1/2}Dg(x)|^4}.$$

Therefore

$$\int_H \varphi'(g(x))\mu(dx)$$

$$= \sum_{k=1}^\infty \lambda_k \int_H D_k\varphi(g(x))\frac{D_k g(x)}{|Q^{1/2}Dg(x)|^2}\mu(dx)$$

$$= -\sum_{k=1}^\infty \lambda_k \int_H \varphi(g(x))\frac{D_k^2 g(x)}{|Q^{1/2}Dg(x)|^2}\mu(dx)$$

$$= -2\sum_{k=1}^\infty \lambda_k \int_H \varphi(g(x))D_k g(x)\frac{\sum_{j=1}^\infty \lambda_j D_k D_j g(x)D_j g(x)}{|Q^{1/2}Dg(x)|^4}\mu(dx)$$

$$+ \sum_{k=1}^\infty \int_H x_k\varphi(g(x))\frac{D_k g(x)}{|Q^{1/2}Dg(x)|^2}\mu(dx).$$

So, (A.5) follows.  □

The following result can be proved similarly.

COROLLARY A.4.  *Assume that Hypothesis A.1 is fulfilled and let $f$ be bounded and Borel. Then $\mu_{fg} \ll \ell$.*

**A.2. Surface measure.** We denote by $K$ the closed set $K = \{g(x) \le 1\}$ and set

$$\Sigma_r = \{g(x) = r\}, \qquad \Sigma = \Sigma_1.$$

We recall the disintegration formula, see, for example, [19, 20]. For any $\varphi : H \to \mathbb{R}$ bounded and Borel we have.

$$(A.7) \qquad \int_H \varphi(x)\mu(dx) = \int_0^{+\infty} \left[ \int_{\Sigma_r} \varphi(x) m_r(dx) \right] \mu_g(dr),$$

where $(m_r)_{r \ge 0}$ is a family of Borel measures on $[0, +\infty)$ such that the support of $m_r$ is included on $\Sigma_r$.

Set

$$\alpha(r) = \int_{\{g \le r\}} d\mu = \mu_g([0, r]).$$

By Proposition A.3 $\alpha$ is a.e. differentiable on $(0, \infty)$. We set

$$\sigma_\mu(\Sigma_r) := \alpha'(r) = \lim_{h \to 0} \frac{1}{2h} \int_{r-h \le g(x) \le r+h} \mu(dx).$$

Now let $f$ bounded and Borel and set

$$\alpha_f(r) = \int_{\{g \le r\}} f \, d\mu = (f\mu)_g([0, r]).$$

Then by Corollary A.4 it follows that $\alpha_f$ is a.e. differentiable. We set

$$\int_{\Sigma_r} f(y)\sigma_{\mu_r}(dy) := \alpha_f'(r) = \lim_{h \to 0} \frac{1}{2h} \int_{r-h \le g(x) \le r+h} f(x)\mu(dx), \qquad \text{a.e. } r > 0.$$

We finally prove.

THEOREM A.5. *Let $f \in B_b(H)$. Then we have*

$$(A.8) \qquad \int_H f(x)\mu(dx) = \int_0^{+\infty} \left[ \int_{\Sigma_r} f(\sigma)\sigma_{\mu_r}(d\sigma) \right] dr.$$

PROOF. Using the disintegration formula (A.7) we have a.e. on $(0, \infty)$

$$\int_{\Sigma_r} f(\sigma)\sigma_{\mu_r}(d\sigma) =: \lim_{h \to 0} \frac{1}{2h} \int_{r-h \le g(x) \le r+h} f(x)\mu(dx)$$

$$= \lim_{h \to 0} \frac{1}{2h} \int_{r-h}^{r+h} \left[ \int_{g^{-1}(r)} f(x) m_r(dx) \right] \sigma_{\mu_r}(\Sigma_r) \, dr.$$

By Lebesgue's theorem we deduce that

$$\int_{\Sigma_r} f(\sigma)\sigma_{\mu_r}(d\sigma) = \int_{g^{-1}(r)} f(x) m_r(dx)\sigma_{\mu_r}(\Sigma_r), \qquad \text{a.e. } r > 0,$$

which yields

$$\int_{g^{-1}(r)} f(x) m_r(dx) = \frac{1}{\sigma_{\mu_r}(\Sigma_r)} \int_{\Sigma_r} f(\sigma) \sigma_{\mu_r}(d\sigma), \qquad \text{a.e. } r > 0.$$

Now the conclusion follows by substituting this into (A.7). $\quad\square$

**Acknowledgment.** We thank the anonymous referee for pertinent observations which have improved the presentation of this work.

V. Barbu
University Al. I. Cuza
Institute of Mathematics Octav Mayer
Blvd. Carol, 11
700506 Iasi
Romania
E-mail: vb41@uaic.ro

G. Da Prato
Scuola Normale Superiore
Piazza dei Cavalieri, 7
56126 Pisa
Italy
E-mail: daprato@sns.it

L. Tubaro
Dipartimento di Matematica
Via Sommarive, 14
38050 Povo (Trento)
Italy
E-mail: tubaro@science.unitn.it